\documentclass[a4paper,12pt]{amsart}

\usepackage{amsthm}

\usepackage{amssymb} 
\usepackage{latexsym} 
\usepackage{amsfonts} 
\usepackage{amsmath} 
\usepackage{eucal} 
\usepackage{bm} 
\usepackage{bbm} 
\usepackage{graphicx} 
\usepackage[english]{varioref} 
\usepackage[nice]{nicefrac} 
\usepackage[all]{xy}

\newcommand{\pos}{\text {\rm pos}}

\newcommand{\ti}{\tilde }

\newcommand{\Ad}{\text {\rm Ad}}

\newcommand{\boa}{\mathbf a}

\newcommand{\kk}{\mathbf k}

\newcommand{\cb}{\mathcal B}

\newcommand{\cp}{\mathcal P}

\newcommand{\cz}{\mathcal Z}

\newcommand{\tG}{\ti G}

\newcommand{\tP}{\ti P}

\def\le{\leqslant}

\def\d{\delta}
\def\D{\Delta}
\def\L{\Lambda}

\def\o{\omega}

\def\s{\sigma}

\def\i{^{-1}}

\theoremstyle{plain}
\newtheorem*{thm*}{Theorem} 

 \newtheorem*{rmk}{Remark}

\newtheorem*{th1}{Lemma 1.1}
\newtheorem*{th2}{Proposition 1.2}
\newtheorem*{th3}{Proposition 1.5}
\newtheorem*{th4}{Lemma 1.7}
\newtheorem*{th5}{Lemma 1.11}
\newtheorem*{th6}{Proposition 1.12}
\newtheorem*{th7}{Lemma 1.13}
\newtheorem*{th8}{Corollary 1.14}
\title{$G$-stable pieces and Lusztig's dimension estimates}
\author{Xuhua He}
\thanks{The author is partially supported by NSF grant DMS-0700589}

\begin{document}

\begin{abstract}
We use $G$-stable pieces to construct some equidimensional varieties and as a consequence, obtain Lusztig's dimension estimates \cite[section 4]{L2}. This is a generalization of \cite{HL}.
\end{abstract}
\maketitle


In Lemma 1.1 and Proposition 1.2, we assume that $G$ is arbitrary connected algebraic group and $\tG$ is an algebraic group with identity component $G$.

\begin{th1} Let $g \in \tG$. Define $i: \tG \to \tG$ by $i(h)=g h g \i h \i$. For any closed subgroup $A$ of $Z_G$ with $g A g \i=A$, set $H_A=\{h \in G; i(h) \in A\}$. Then

(1) $H_A$ is an algebraic group and $i|_{H_A}: H_A \to A$ is a morphism of algebraic groups.

(2) $i(A)^0=i(H_A)^0$.

(3) $\dim(H_A)=\dim(Z_G(g))+\dim(A)-\dim(Z_A(g))$.
\end{th1}

If $h, h' \in H_A$, then \begin{align*} i(h h') &=g h h' g \i (h') \i h \i=(g h g \i h \i) h (g h' g \i (h') \i) h \i \\ &=i(h) h i(h') h \i=i(h) i(h') \in A \end{align*} 
and $h h' \in H_A$. If $h \in H_A$, then $i(h \i)=h \i i(h) \i h=i(h) \i \in A$ and $h \i \in H_A$. Part (1) is proved.

Now $i(A)^0$ is a connected subgroup of $i(H_A)$. Define $\d: A \to A$ by $\d(z)=g z g \i$. Then $$\dim(i(A))=\dim(A)-\dim(A^\d).$$ Define $\s: A \to A$ by $\s(z)=\d^{m-1}(z) \d^{m-2}(z) \cdots z$, where $m$ is the order of the automorphism $\d$. Then $\s$ is a group homomorphism and $$i(H_A) \subset \{z \in Z; \s(z)=1\}.$$ Notice that $\s(A^\d)=\{t^m; t \in A^\d\}$ is of dimension $\dim(A^\d)$. Thus 
\begin{align*}
\dim(i(H_A)) & \le \dim(A)-\dim(\s(A)) \le \dim(A)-\dim(\s(A^\d)) \\ &=\dim(A)-\dim(A^\d).
\end{align*}

Therefore, $\dim(i(A))=\dim(i(H_A))=\dim(A)-\dim(A^\d)$. Part (2) is proved.

Since the kernal of $i|_{H_A}$ is $Z_G(g)$, \begin{align*} \dim(H_A) &=\dim(Z_G(g))+\dim(i(H_A)) \\ &=\dim(Z_G(g))+\dim(A)-\dim(A^\d).\end{align*} Part (3) is proved. \qed

\begin{th2} Let $D, D'$ be connected group of $\tG$. Define $\d: Z_G \to Z_G$ by $\d(z)=g z g \i$ for any $g \in D$ and $\d': Z_G \to Z_G$ by $\d'(z)=g' z (g') \i$ for any $g' \in D'$. Let $c$ be a $G$-conjugacy class in $D$ and $Z$ be a closed subgroup of $Z_G$ with $\d(Z)=Z$. Set $$X=\{(g, g'); g \in c Z, g' \in D', g g' g \i (g') \i \in Z\}.$$ If $X \neq \emptyset$, then $D D'=D' D$ and $\d \d'(z)=\d' \d(z)$ for $z \in Z_G$. Moreover, $X$ is of pure dimension $\dim(G)+\dim(Z)-\dim(\frac{\d'(Z)}{Z \cap \d'(Z)})^\d$.
\end{th2}

It is easy to see that if $X \neq \emptyset$, then $D D'=D' D$. Thus for $g \in D$ and $g' \in D'$, $g \i (g') \i g g' \in G$ and $\d \i (\d') \i \d \d'(z)=z$ for all $z \in Z_G$.

Consider the projection map $X \to D$ defined by $(g, g') \mapsto g$. Let $g$ be in the image and $X_g$ be the fiber over $g$. Fix $g' \in X_g$. Set $i: \tG \to \tG$ by $i(h)=g h g \i h \i$. Then $i(h g')=i(h) i(g')$ for $h \in G$. Hence $X_g=H_Z g'$.

Let $z \in Z$. Then for $h \in G$, $$(z g) (h g') (z g) \i (h g') \i=z i(h) i(g') \d'(z) \i.$$ If $z i(h) i(g') \d'(z) \i \in Z$, then $i(h) \in Z \d'(Z)$. Hence $h \in H_{Z \d'(Z)}$ and $\d'(z) \in i(H_{Z \d'(Z)}) Z$. On the other hand, if $z \in Z$ with $\d'(z) \in i(H_{Z \d'(Z)}) Z$, then there exists $h \in G$ with $z i(h) i(g') \d'(z) \i \in Z$. Therefore, for $z \in Z$, $X_{z g} \neq \emptyset$ if and only if $z \in Z'$, where $Z'=\{z \in Z; \d'(z) \in i(H_{Z \d'(Z)}) Z\}$. It is easy to see that $Z'$ is an algebraic group.

By part (2) of the previous lemma, $i(H_{Z \d'(Z)})^0 \subset i(Z \d'(Z))$. Hence $\d'(Z')^0=\bigl(\d'(Z) \cap i(Z \d'(Z)) Z \bigr)^0=\bigl(\d'(Z) \cap i(\d'(Z)) Z \bigr)^0$. Notice that $i(Z) \subset \d(Z) Z=Z$ and $i \d'(Z) \subset \d \d'(Z) \d'(Z)=\d'(Z)$. Now $i$ induces a group morphism $\bar i: \d'(Z)/(Z \cap \d'(Z)) \to \d'(Z)/(Z \cap \d'(Z))$. We have that $z \in \d'(Z) \cap i(\d'(Z)) Z$ if and only if $z (Z \cap \d'(Z))$ is contained in the image of $\bar i$. Hence $\dim(Z')=\dim(\frac{\d'(Z)}{Z \cap \d'(Z)})-\dim(\frac{\d'(Z)}{Z \cap \d'(Z)})^\d+\dim(Z \cap \d'(Z))=\dim(Z)-\dim(\frac{\d'(Z)}{Z \cap \d'(Z)})^\d$.

Set $Y=\{(h, z); z \in Z', h \in X_{z g}\}$. Then we have the projection map $Y \to Z'$ and each fiber is isomorphic to $H_Z$. Hence $Y$ is of pure dimension $\dim(Z')+\dim(H_Z)$.

Consider the morphism $G \times Y \to X$ defined by $$(g_1, h, z) \mapsto (g_1 h g_1 \i, z g_1 g g_1 \i).$$ Then this morphism is surjective and the fiber over $(h, z g)$ is $$\{(g_1, g_1 \i h g_1, i(g_1) z); g_1 \in H_Z\},$$ which is of pure dimension $\dim(H_Z)$. Therefore $X$ is of pure dimension $\dim(G)+\dim(Z')$. \qed

\subsection*{1.3} From now on, we assume that $G$ be a simply-connected, semisimple algebraic group over an algebraically closed field $\kk$. We fix a Borel subgroup $B$ of $G$ and a maximal torus $T \subset B$. Let $I$ be the set of simple roots determined by $B$ and $T$.

For any $J \subset I$, let $P_J$ be the standard parabolic subgroup corresponding to $J$ and $\cp_J$ be the set of parabolic subgroups that are $G$-conjugate to $P_J$. We simply write $\cp_{\emptyset}$ as $\cb$. Let $L_J$ be the Levi subgroup of $P_J$ that contains $T$.

For any parabolic subgroup $P$, let $U_P$ be the unipotent radical of $P$ and $H_P$ be the inverse image of the connected center of $P/U_P$ under the projection map $\pi_P: P \to P/U_P$. We simply write $U$ for $U_B$.

For $J \subset I$, we denote by $W_J$ the standard parabolic subgroup of $W$ generated by $J$ and by $W^J$ (resp. ${}^J W$) the set of minimal coset representatives in $W/W_J$ (resp. $W_J \backslash W$). For $J, K \subset I$, we simply write $W^J \cap {}^K W$ as ${}^K W^J$.

For $P \in \cp_J$ and $Q \in \cp_K$, we write $\pos(P, Q)=w$ if $w \in {}^J W^K$ and there exists $g \in G$ such that $P=g P_J g \i$, $Q=g \dot w P_K \dot w \i g \i$, where $\dot w$ is a representative of $w$ in $N(T)$.

For $g \in G$ and $H \subset G$, we write $^g H$ for $g H g \i$.

For any algebraic group $H$, let $H^0$ be its identity component.

\subsection*{1.4} Let $\s$ be a diagram automorphism of $G$, i.e., an automorphism of $G$ that stabilizes $B$ and $T$ and the order of $\s$ as an automorphism of $G$ coincides with the order of the associated permutation on $I$. We use the same symbol $\s$ for the associated automorphism on $W$ and associated permutation on $I$. Set $\tG=G \rtimes <\s>$, where $<\s>$ is the finite subgroup of $G$ generated by $\s$. We simply write the element $(g, \s^n) \in \tG$ as $g \s^n$. For each element $g \in \tG$, we write $g_s$ for its semisimple part and $g_u$ its unipotent part.

Let $D=(G, \s)$ be a connected component of $\tG$. We have the following result.

\begin{th3} Let $g \in D$. Then $g$ is $G^0$-conjugate to an element of the form $t \s u$, where $t \in (T^\s)^0$ and $u$ is a unipotent element in $Z_G^0(t \s_s)$.
\end{th3}

By \cite[Lemma 7.3]{St}, after conjugate by an element in $G^0$, we may assume that $g \in B \s \subset B \ltimes <\s>$. Then $g_s \in B \s_s$ and $g_u \in B \s_u$. Then after conjugate by an element in $B$, we may assume that $g_s=t_1 \s_s$ and $g_u=t_2 \s_u u$ for $t_1, t_2 \in T$ and $u \in U$. By \cite[1.2]{L2}, after conjugate by an element in $T$, we may assume furthermore that $t_2 \in (T^{\s_u})^0$ . Consider the group morphism $B \ltimes <\s> \rightarrow T \ltimes <\s>$. Since $g_u$ is unipotent, then $t_2 \s_u$ is also unipotent. Notice that $t_2$ commutes with $\s_u$. Then $t_2$ is unipotent and $t_2=1$.

Since $\s$ is a diagram automorphism, $\s_s$ and $\s_u$ are also diagram automorphisms. In particular, $\rho^\vee(t) \in T^{\s_\s} \cap T^{\s_u}$ for all $t \in \kk^\times$. Hence $\Ad(\rho^\vee(t)) g_u \in Z_G(g_s)$. Since $\s_u$ is contained in the closure of $\{\Ad(\rho^\vee(t)) g_u)\}$ and $Z_G(g_s)$ is closed, we have that $\s_u \in Z_G(g_s)$. We also have that $u \in Z^0_G(g_s)$.

Now $\s_u$ commutes with $t_1 \s_s$. Hence $t_1 \in T^{\s_u}$ By \cite[9.6]{Bo}, $(T^{\s_u})$ is connected. Notice that $\s_s$ is an automorphism on $T^{\s_u}$. Then by \cite[1.2]{L2}, after conjugate by $T^{\s_u}$, $t_1 \in ((T^{\s_u})^{\s_s})^0 \subset (T^\s)^0$. \qed

\subsection*{1.6}\label{1} Let $D// G$ be set of closed $G$-conjugacy classes in $D$. By geometric invariant theory, $D//G$ has a natural structure of affine variety and there is a well-defined morphism ${\rm St}: D \to D//G$ which maps the element $g \in D$ to the unique closed $G$-conjugacy class in $D$ that is contained in the closure of the $G$-conjugacy class of $g$. If $\s$ is trivial, then ${\rm St}$ is just the Steinberg morphism of $G$. Thus for arbitrary $\s$, we call ${\rm St}$ the Steinberg morphism of $D$ and the fibers the Steinberg fibers of $D$.

By the previous proposition, any element $g \in D$ is of the form $t \s u$, where $t \in (T^\s)^0$ and $u$ is a unipotent element in $Z^0(t \s_s)$. Moreover, $t \s_s$ is contained in the closure of the $G$-conjugacy class of $g$. Hence ${\rm St}(g)={\rm St}(t \s_s)$. Notice that $t \sigma$ is quasi-semisimple in the sense of \cite[Sect.9]{St}, i.e. the automorphism of $G$ obtained by conjugation by $t \s_s$ will fix a Borel subgroup and a maximal torus thereof. As a consequence, the $G$-conjugacy class of $t \s_s$ in $D$ is closed \cite[II.1.15(f)]{Sp}. We conclude that any Steinberg fiber of $D$ is of the form
$$\bigcup_{g \in G} \bigcup_{u \text{ is unipotent in } Z_G(t \s_s)^0} g (t \s u) g \i$$ for some $t \in (T^{\s})^0$. It is known that $Z_G(t \s_s)^0$ is reductive and the set of unipotent elements in a reductive group is an irreducible variety. Thus 

(a) each Steinberg fiber is irreducible. 

Moreover, there are only finitely many unipotent conjugacy classes in a reductive group \cite{L1}. Therefore 

(b) each Steinberg fiber contains finitely many $G$-conjugacy classes.

\begin{th4} Let $\boa$ be a Steinberg fiber in $D$ and $J \subset I$ with $\s(J)=J$.
Then there exists finitely many $L_J$-conjugacy classes $c_1, \cdots, c_m$ in $N_{\tG}(P_J) \cap N_{\tG}(L_J) \cap D$ such that $\{g; g \in N_{\tG}(P_J) \cap D, St(g)=\boa\}=\sqcup_i c_i U_{P_J}$.
\end{th4}

Let $l \in N_{\tG}(P_J) \cap N_{\tG}(L_J)$ and $u \in U_{P_J}$. Then it is easy to see that $l$ is contained in the closure of $\{t l u t \i; t \in (Z(L_J)^\s)^0\}$. Hence $St(l)=St(l u)$. In other words, \begin{align*} & \{g; g \in N_{\tG}(P_J) \cap D, St(g)=\boa\} \\ &=\{l u; l \in N_{\tG}(P_J) \cap N_{\tG}(L_J) \cap D, u \in U_{P_J}, St(l)=\boa\}.\end{align*}

By \cite[Proposition 1.14]{L2}, any quasi-semisimple element in $D$ (resp. $N_{\tG} (P_J) \cap N_{\tG} (L_J) \cap D$) is $G$-conjugate (resp. $L_J$-conjugate) to $T_1 \s$, where $T_1=(T^\s)^0$. Notice that $\{t \in T_1; St(t \s)=\boa\}$ is a finite set. Hence there are only finitely many quasi-semisimple $L_J$-conjugacy classes in $N_{\tG} (P_J) \cap N_{\tG} (L_J) \cap D$ that are contained in $St \i(\boa)$. One can see that a $L_J$-conjugacy class $c$ is contained in $St \i (\boa)$ if and only if the unique $L_J$-conjugacy class $c'$ that is contained in $c$ is also contained in $St \i (\boa)$. Then the lemma follows from Lemma 1.7 (b). \qed

\subsection*{1.8}\label{3} For $J \subset I$, set \begin{align*} Z_J &=\{(P, P', g U_P); P, P' \in \cp_J, g \in G, P'={}^g P\}, \\ Z'_J &=\{(P, P', g H_P); P, P'
\in \cp_J, g \in G, P'={}^g P\}
\end{align*} with the $G \times G$-action defined by \begin{align*} (g_1, g_2) \cdot (P, P', g U_P) &=({}^{g_2} P, {}^{g_1} P', g_1 g U_P g_2 \i), \\ (g_1, g_2) \cdot (P, P', g H_P) &=({}^{g_2} P, {}^{g_1} P', g_1 g H_P g_2 \i). \end{align*}

Set $h_J=(P_J, P_J, U_{P_J}) \in Z_J$ and $h'_J=(P_J, P_J, H_{P_J}) \in Z'_J$. By \cite[section 3]{L3}, \cite[section 1]{H1} and the remark of \cite[Corollary 5.4]{H2}, we have partitions \[\tag{a} Z_J=\sqcup_{w \in {}^J W} Z_{J; w} \qquad \text{ and } \qquad Z'_J=\sqcup_{w \in {}^J W} Z'_{J; w}, \] where $Z_{J; w}=G_{\D} \cdot (B w B, B) h_J$ and $Z'_{J; w}=G_{\D} \cdot (B w B, B) h'_J$. The subvarieties $Z_{J; w}$ (resp. $Z'_{J; w}$) are called $G$-stable pieces of $Z_J$ (resp. $Z'_J$).

Fix $w \in {}^J W$. Let $K=I(J, id; w)$. Then by \cite[3.14]{L3},\\
(b) there is a canonical bijection between the $G_{\D}$-orbits on $Z_{J, w}$ and the $L_K$-conjugacy classes on $w L_K$.

By \cite[section 3]{L3}, we have $G$-equivariant morphisms $pr: Z_{J; w} \to G/P_K$ and $pr': Z'_{J; w} \to G/P_K$, where $G$ acts diagonally  on $Z_{J; w}$ and $Z'_{J; w}$ and acts in the natural way on $G/P_K$. Moreover, by \cite[Proposition 1.10]{H1}, \begin{align*} \tag{c} & pr(z)=P_K \text{ if and only if }z=(p w, 1) \cdot h_J \text{ for some } p \in P_K, \\ & pr'(z)=P_K \text{ if and only if } z=(p w, 1) \cdot h'_J \text{ for some } p \in P_K.\end{align*}

Also we have that the closure of any $G$-stable piece is a union of $G$-stable pieces.

(d) $\overline{Z_{J; w}}=\sqcup_{w' \in {}^J, w' \le_{J, id} w} Z_{J, w'}$. See \cite[Proposition 4.6]{H1} and \cite[Proposition 5.8]{H2}.

\subsection*{1.9} If $\s(J)=J$, then the action of $G \times G$ on $Z_J$ and $Z'_J$ can be extended in a natural way to an action of $\tG \times \tG$.

Now set \begin{gather*} \L_{J, D}=\{(z, g) \in Z_J \times D; (g, g) \cdot z=z\}, \\ \L'_{J, D}=\{(z, g) \in Z'_J \times D; (g, g) \cdot z=z\}.
\end{gather*} For $w \in W^J$, set $\L_{J, D; w}=\{(z, g) \in \L_{J, D}; z \in Z_{J; w}\}$ and $\L'_{J, D; w}=\{(z, g) \in \L'_{J, D}; z \in Z'_{J; w}\}$.

\subsection*{1.10} Set $\tP=N_{\tG} P$ for any parabolic subgroup $P$ of $G$. Define the action of $P_J$ on $G \times \tP_J$ by $p \cdot (g, p')=(g p \i, p p' p \i)$. Let $G \times_{P_J} \tP_J$ be the quotient space. Then we may identify $G \times_{P_J} \tP_J$ with $\{(P, g); P
\in \cp_J, g \in \tP\}$ via $(g, p) \mapsto ({}^g P_J, g p g \i)$.

Let $c$ be a subvariety of $N_{\tG}(P_J) \cap N_{\tG}(L_J) \cap D$
that is stable under the conjugation action of $L_J$. Then $c
U_{P_J}$ and $c H_{P_J}$ are stable under the conjugation action of
$P_J$. So we may define $G \times_{P_J} c U_{P_J} \subset G
\times_{P_J} c H_{P_J} \subset G \times_{P_J} \tP_J$.

Now set \begin{align*} \L_{J, c} &=\{(P, P', g U_P, h) \in \L_{J,
D}; (P, h), (P', h) \in G \times_{P_J} c U_{P_J}\}, \\
\L'_{J, c} &=\{(P, P', g H_P, h) \in \L'_{J, D}; (P, h), (P', h) \in G
\times_{P_J} c H_{P_J}\}.
\end{align*}
For $w \in W^J$, set $\L_{J, c; w}=\L_{J, D; w} \cap \L_{J, c}$ and $\L'_{J, c; w}=\L'_{J, D; w} \cap \L'_{J, c}$.

\begin{th5} Let $w \in {}^J W$ and $K=I(J, id; w)$. Let $c$ be a $L_K$-conjugacy class in $N_{\tG}(P_K) \cap N_{\tG}(L_K) \cap D$.
Set $X_{J, c; w}=\{(z, g) \in \L_{J, c; w}; pr(z)=P_K\}$ and
$X'_{J, c; w}=\{(z, g) \in \L'_{J, c; w}; pr'(z)=P_K\}$. Then 

(1) If $X_{J, c; w} \neq \emptyset$, then $X_{J, c; w}$ is of
pure dimension $\dim(P_K)$.

(2) If $X'_{J, c; w} \neq \emptyset$, then $X'_{J, c; w}$ is of pure dimension $$\dim(P_K)-\dim(\frac{\Ad(w) Z^0(L_J)}{Z^0(L_J) \cap \Ad(w) Z^0(L_J)})^\s.$$
\end{th5}

\begin{rmk} (1) If $X_{J, c; w} \neq \emptyset$ or $X'_{J, c; w} \neq \emptyset$, then as we will see in the proof, $\s \Ad(w) (z)=\Ad(w) \s(z)$ for $z \in Z^0(L_J)$. So we have that  $\s Z^0(L_J)=Z^0(L_J)$ and $\s \Ad(w) Z^0(L_J)=\Ad(w) Z^0(L_J)$. Therefore $(\frac{\Ad(w) Z^0(L_J)}{Z^0(L_J) \cap \Ad(w) Z^0(L_J)})^\s$ is defined.

(2) If $w (J)=J$, then $\Ad(w) Z^0(L_J)=Z^0(L_J)$ and $$\dim(\frac{\Ad(w) Z^0(L_J)}{Z^0(L_J) \cap \Ad(w) Z^0(L_J)})^\s=0.$$ 

If $w (J) \neq J$, then there exists $j \in J$ such that $w (j)$ is not spanned by the simple roots in $J$. Set $z_t=\prod_{j \notin J} \o_j(t)$ for $t \in \kk^*$. Then $z_t \notin \Ad(w) \i Z^0(L_J)$ for all $t \in \kk^*-\{1\}$. Moreover, $\s \Ad(w) (z_t)=\Ad(w) \s(z_t)=\Ad(w)(z_t)$. Hence $\Ad(w) z_t \in (\Ad(w) Z^0(L_J))^\s$.

So $\dim(\frac{\Ad(w) Z^0(L_J)}{Z^0(L_J) \cap \Ad(w) Z^0(L_J)})^\s=0$ if and only if $w (J)=J$.
\end{rmk}

We only prove part (1) here. Part (2) can be proved in a similar way.

Let $\bar p: P_K \to L_K$ be the projection map. By \ref{3} (c), $pr(z)=P_K$ if and only if $z=(p w, 1) \cdot h_J$ for some $p \in P_K$. Moreover, the morphism $f: \{z \in Z_{J; w}; pr(z)=P_K\} \to L_K$ defined by $(p w, 1) \cdot h_{J, D} \mapsto \bar p(p)$ is well-defined. Now consider the morphism $X_{J, c; w} \to L_K \times c$ which sends $(z, g)$ to $(f(z), \bar p(g))$. We see that the image is contained in $\{(l, l') \in L_K \times c; l w l'=l' l w\}$, which is of pure dimension $\dim(L_K)$.

Let $l \in L_K$ and $l' \in c$ with $l w l'=l' l w$. Then the fiber $Y$ over $(l, l')$ is $\{(z, u l'); z \in (U_{P_K} l w, 1) \cdot h_J, u \in U_{P_K}, (u l', u l') \cdot z=z\}$. Define the action of $U_{P_K}$ on $Y$ by $u_1 \cdot (z, u l')=((u_1, u_1) \cdot z, u_1 u l' u_1 \i)$. Then the projection map $Y \to (U_{P_K} l w, 1) \cdot h_J=(U_{P_K})_{\D} (l
w, 1) \cdot h_J$ is $U_{P_K}$-equivariant for the diagonal $U_{P_K}$-action on $(U_{P_K} l w, 1) \cdot h_J$. Since $U_{P_K}$ acts transitively on $(U_{P_K} l w, 1) \cdot h_J$, the projection map is a locally trivialy fibration with fibers isomorphic to $\{u \in U_{P_K}; (u
l w, u) \cdot h_J=(l w, 1) \cdot h_J\}$. In particular, $Y$ is of pure dimension $\dim(U_{P_K})$. Therefore $X_{J, c; w}$ is of pure dimenion $\dim(L_K)+\dim(U_{P_K})=\dim(P_K)$. \qed

\begin{th6} Let $\boa$ be a Steinberg fiber of $N_{\tG}(P_J) \cap N_{\tG}(L_J) \cap D$.

(1) If $\L_{J, \boa; w} \neq \emptyset$, then $\L_{J, \boa; w}$ is of pure dimension $\dim(G)$.

(2) If $\L'_{J, \boa; w} \neq \emptyset$, then $\L'_{J, \boa; w}$ is of pure dimension $$\dim(G)-\dim(\frac{\Ad(w) Z^0(L_J)}{Z^0(L_J) \cap \Ad(w) Z^0(L_J)})^\s.$$
\end{th6}

We only prove part (1) here. Part (2) can be proved in a similar
way.

Let $K=I(J, id; w)$. For $(z, g) \in X_{J, \boa; w}$, we have that $^g P_K=g \cdot pr(z)=pr\bigl( (g, g) \cdot z \bigr)=pr(z)=P_K$. Hence $g \in P_K$. By lemma 1.7, $X_{J, \boa; w}=\sqcup_i X_{J, c_i; w}$ for finitely many $L_K$-conjugacy classes $c_i$. By lemma 1.11, $X_{J, \boa; w}$ is of pure dimension $\dim(P_K)$.

Let $\pi: \L_{J, a; w} \to Z_{J; w}$ be the projection map. It is easy to see that $\pi$ is $G$-equivariant for the diagonal $G$-action. Thus $pr \circ \pi: \L_{J, \boa; w} \to G/P_K$ is also $G$-equivariant. Since $G$ acts transitively on $G/P_K$, $pr \circ \pi$ is a locally trivial fibration with fibers isomorphism to $X_{J, \boa; w}$. Thus $\L_{J, w; \boa}$
is of pure dimension $\dim(G)$. \qed

\begin{th7} Let $c$ be a $L_J$-conjugacy class in $N_{\tG}(P_J) \cap N_{\tG}(L_J) \cap D$.
Set \begin{align*} \cz_{J, c} &=\{(P, P', g); (P, g), (P', g) \in G
\times_{P_J} c U_{P_J}\}, \\
\cz'_{J, c} &=\{(P, P', g); (P, g), (P', g) \in G \times_{P_J} c H_{P_J}\}.\end{align*}

(1) Define the map $\L_{J, c} \to \cz_{J, c}$ by $(P, P', k U_P, g) \mapsto (P, P', g)$. If $\cz_{J, c} \neq \emptyset$, then the map is surjective and each fiber is of pure dimension $\dim(L_J)-\dim(c)$.

(2) Define the map $\L'_{J, c} \to \cz'_{J, c}$ by $(P, P', k H_P, g) \mapsto (P, P', g)$. If $\cz'_{J, c} \neq \emptyset$, then the map is surjective and each fiber is of pure dimension $\dim(L_J)-\dim(c)-\dim(Z(L_J)^\s)$.
\end{th7}

We only prove part (1) here. Part (2) can be proved in the same way.

Let $(P, P', g) \in \cz_{J, c}$. Then there exists $k \in G$ such that $P'={}^k P$. By definition, $k g k \i U_{P'}$ and $g U_{P'}$ are $P'$-conjugate. Therefore, there exists $l \in P'$, such that $l k g k \i l \i \in g U_{P'}$. In other words, $(g, g) \cdot (P, P', l
k U_P)=(P, P', l k U_P)$. So the map is surjective.

Assume that $(P, P', k U_P, g), (P, P', k' U_P, g) \in \L_{J, c}$. Then $k \i k' \in P$ and $(k \i k') g U_P (k \i k') \i=g U_P$. Thus the fibers of the map $\L_{J, c} \to \cz_{J, c}$ are isomorphic to $\{(l U_P; l \in P; l g U_P l \i=g U_P\}$. Notice that $(P, g) \in G \times_{P_J} c U_{P_J}$.  is of dimension $\dim(L_J)-\dim(c)$. \qed

\

Now combining Proposition 1.12 and Lemma 1.13, we have the following result.

\begin{th8} Let $c$ be a $L_J$-conjugacy class in $N_{\tG}(P_J) \cap N_{\tG}(L_J) \cap D'$.
Then

(1) $\dim(\cz_{J, c}) \le \dim(G)-\dim(L_J)+\dim(c)$.

(2) $\dim(\cz'_{J, c}) \le \dim(G)-\dim(L_J)+\dim(c)+\dim(Z(L_J)^\s)$. More
precisely, define $\cz'_{J, c; w}=\{(P, P', g) \in \cz'_{J, c}; (P, P') \in G_{\D} \cdot (P_J, {}^w P_{J'})\}$ for $w \in {}^J W^{J'}$. Then \begin{align*} \dim(\cz'_{J, c; w}) \le
\dim(G) & -\dim(L_J)+\dim(c)+\dim(Z(L_J)^\s) \\ &-\dim(\frac{\Ad(w) Z^0(L_J)}{Z^0(L_J) \cap \Ad(w) Z^0(L_J)})^\s.\end{align*}
\end{th8}

\begin{rmk}
Part (1) was first proved in \cite[Proposition 4.2 (d)]{L2}. By the remark of Lemma 1.11, part (2) is a stronger version of \cite[Proposition 4.2(c)]{L2}.
\end{rmk}

\bibliographystyle{amsalpha}

\end{document}